\documentclass[12pt]{amsart}
\usepackage{amssymb}
\usepackage{amsfonts}
\usepackage{amssymb}
\usepackage{amsthm}

\usepackage{amsmath}

\usepackage{calc}

\usepackage{times}

\usepackage{amscd}
\theoremstyle{plain}
\newtheorem{thm}{Theorem}
\newtheorem*{thma}{Theorem A}
\newtheorem*{thmb}{Theorem B}
\newtheorem*{thmno}{Theorem}
\newtheorem{lem}[thm]{Lemma}
\newtheorem{prop}[thm]{Proposition}
\newtheorem{cor}[thm]{Corollary}

\theoremstyle{definition}
 
\newtheorem*{remno}{Remarks}

\def\PP{{\mathbb P}}
\def\co{{\mathcal O}}
\def\ce{{\mathcal E}}

\def\Pn{\PP ^n}
\def\deg{\mathop{\rm deg}}
\def\cl{{\mathcal L}}
\def\ind{\indent\indent}
\def\rank{\mathop{\rm rank}}
\def\bigskip{\vskip24pt}
\def\medskip{\vskip12pt}
\def\Pic{\mathop{\rm Pic}}

\title{ON MANIFOLDS OF SMALL DEGREE} 
\author{PALTIN IONESCU}
\begin{document}
\def\leq{\leqslant}
\def\geq{\geqslant}
\def\eq{\leqslant}
\def\eq{\geqslant}
\address{University of Bucharest, Department of Mathematics, 14 Academiei Street, 70109 Bucharest, Romania}
\address{and} 
\address{Institute of Mathematics ``Simion Stoilow" of the Romanian Academy, P.O. Box 1--764, 014700 Bucharest, Romania; }
\address{E-mail: {\rm Paltin.Ionescu@imar.ro}}

\subjclass{Primary: 14N25, 14N30; Secondary: 14M99, 14J45.}
\begin{abstract}
Let $X\subset \Pn $ be a complex projective manifold of degree
$d$ and arbitrary dimension. The main result of this paper gives a
classification of such manifolds (assumed moreover to be connected,
non-degenerate and linearly normal) in case $d\leq n$. As a by-product of
the classification it follows that these manifolds are either rational
or Fano. In particular, they are simply connected (hence regular) and of
negative Kodaira dimension. Moreover, easy examples show that $d\leq
n$ is the best possible bound for such properties to hold true. The proof of
our theorem makes essential use of the adjunction
mapping and, in particular, the main result of \cite{11} plays a crucial role
in the argument.
\end{abstract}
\thispagestyle{empty}

\maketitle
\thispagestyle{empty}

\section{Introduction} \label{S0}

Let $X\subset \PP^n$ be a complex connected projective 
manifold of dimension $r$ and degree $d$. Assume moreover that $X$ is non-degenerate and $d\leq n$. The results contained in this paper have the following topological consequence:
$$\hbox{\it  If $X$ is as above, then $X$ is simply connected}\leqno{(*)}$$
The bound $d\leq n$ is optimal for the validity of $(*)$. Indeed, there exist $r$-dimensional elliptic scrolls in $\PP^{2r}$, of degree $2r+1$ (see \cite{10}, 5.2); they have $b_1=2$.

To the best of our knowledge, $(*)$ was not even conjectured before. There was, however, the following question raised by F.L. Zak:

\smallskip
{\it Is a linearly normal $r$-dimensional manifold in $\PP^{2r+1}$ of degree $\leq 2r+1$ and whose embedded secant variety equals $\PP^{2r+1}$, a regular variety (i.e.\ having $b_1=0$)?}
\smallskip
 
 We refer the interested reader to \cite {1} for a pertinent discussion about
the relevance of Zak's question. It follows from $(*)$ that the answer to this question is positive, even under more general assumptions.
 
 We would like to mention also two related topological  ancestors of $(*)$.
 The first one is  (a special case of)
Barth-Larsen's theorem (see \cite{BL} and, for a singular version,~\cite{FL}):
\begin{equation*}\hbox{\it If $r\geq {\rm codim} _{\PP^n}(X) +1$,   then $\pi_1(X)=(0)$.}\tag{\rm B\hbox{-}L}\end{equation*}

The second result is Gaffney-Fulton-Lazarsfeld's theorem about branched 
coverings of $\PP^r$ (see \cite{GL, FL}):
\begin{equation*}\begin{split}&\mbox{\it If $X\to \PP^r$ is a normal finite covering of degree $d\leq r$,}\\ &\mbox{\it then $\pi_1(X)=(0)$.}\end{split}\tag{\rm G\hbox{-}F\hbox{-}L}\end{equation*}
Note that, for $d\leq r$, $(*)$ follows either from (B-L) or from (G-F-L).
We refer to \cite{FL} for a very nice discussion of such topological aspects.
Recall that the $\Delta$-genus of $X$ in $\PP^n$ is, by definition, the non-negative integer $\Delta:= d+r- h^0(X,\mathcal{O}_X(1))$. Assuming $X$ to be linearly normal in $\PP^n$ (which is not restrictive), condition $d\leq n$ may be restated as:
$$r\geq \Delta +1.$$  
So we see that  $(*)$ is a Barth-Larsen-type result in which codimension is
replaced by $\Delta$-genus.

Our proof of $(*)$ is, however, not topological. We deduce $(*)$ from the following geometric result:
\begin{equation*}\begin{split}&\hbox{\it  If $X$ is as above, then either:\hskip170pt }\\
 &\qquad\hbox{(1) {\it $b_2=1$ and $X$ is a Fano manifold, or}}\\
 &\qquad\hbox{(2) {\it $b_2\geq2$ and $X$ is rational.}}\hskip30pt
 \end{split}\tag{$**$}\end{equation*}
 
 It is well-known that both rational and Fano manifolds are simply-\break connected; see \cite{KMM} for a far-reaching common generalization. So $(*)$ follows from $(**)$. The first case in $(**)$ may be seen as generic, as it includes all complete intersections of dimension at least three. Indeed, we shall prove:
\begin{equation*}\begin{split} &\mbox{\it Manifolds with $d\leq n$ and $b_2 \geq2$ may be classified completely. }\\ &\mbox{\it  There are $6$ infinite series (having arbitrarily large dimension } \\&\mbox{\it  and degree) 
and $14$ ``sporadic" examples.
Moreover, all turn}\\&\mbox{\it    out to be rational.}\end{split}\tag{$*\!*\!*$}\end{equation*}

The precise list is given in the statement of the main result, see the next section.

The proof of the main theorem will occupy Section~\ref{S3}. 
It relies on a very detailed study of the adjunction mapping (see e.g.\ \cite{3}, Chapters 9--11 for a complete treatment). Moreover, the main result of \cite{11} plays a key role in the proof. We note that, besides classical adjunction theory,  some nontrivial facts coming from Mori theory are also used in \cite{11}. Finally, the classification of manifolds of small $\Delta$-genus (cf.\ \cite{5}, \cite{6}, \cite{9}) is also needed.
\smallskip

The present work is a slightly revised version of a paper with the same title that was circulated as Preprint no.\ 17, IMAR, Bucharest, December~2000.   
 
 \section{Statement of the main result} \label{S1}
 
Our main result is the following:

\begin{thmno}\label{th1} Let $X\subset \Pn$ be a connected projective manifold over 
$\mathbb{C}$, of dimension $r$ and degree $d$. Assume moreover that $X$ is 
non-degenerate and linearly  normal. If $d\leq n$, then one of the
following holds:

\hskip6pt{\rm (i)} $r\geq 1$, $X$ is Fano, $b_2(X)=1$;

\hskip3pt{\rm (ii)} $X$ is Fano and either:

\indent\indent {\rm (a)} $2\leq r\leq 4$, $3\leq d\leq8$, $X$ is a
classical del Pezzo manifold with $b_2(X) \geq 2$ {\rm (cf. Theorem B below)}\/;

\indent\indent{\rm (b)} $r=3$, $d=9$, $X$ is the Segre embedding of
$\PP ^1\times \mathbb{F}_1$, where $\mathbb{F}_1$ is the blowing-up of
$\PP ^2$ in a point, embedded in $\PP ^4$ as a rational scroll
of degree $3$;

\indent\indent{\rm (c)} $X$ is one of the following scrolls over 
$\mathbb{P}^2$: $r=4$, $d=10$, $X\simeq \PP  (T_{\PP ^2} 
\oplus \co_{\PP ^2}(1))$ or $r=4$, $d=11$, $X\simeq \PP (\co _{\mathbb{
P}^2}(1)\oplus \co _{\PP ^2}(1)\oplus \co _{\PP ^2}(2))$ or $r=5$,
$d=10$, $X$ is the Segre embedding of $\PP^2\times \PP^3$;

{\rm (iii)} $r\geq 2$, $d\geq r$, $X$ is a scroll over $\PP ^1$ (i.e. a linear 
section of the Segre embedding of $\PP ^1\times \PP^{m}$); 

\hskip.5pt{\rm (iv)} $r\geq 3$ and there is a  vector bundle 
$\ce $ over $\PP ^1$, of 
rank $r+1$ and of splitting type $e=(e_0, \ldots,e_r)$, such that, if 
$L$ 
denotes the tautological divisor on $\PP (\ce )$ and $F$ denotes a fibre 
of the projection $\PP (\ce) \to \PP ^1$, $X$ embeds in $\PP (\ce )$, 
$L|_X=H$ and either:

\ind {\rm (a)} $n=d=2r-1$, $e=(1, \ldots, 1, 0,0)$, $X\in
|2L+F|$;

\ind {\rm (b)} $n=d=2r$, $e=(1,\ldots,1,0)$, $X\in |2L|$;

\ind {\rm (c)} $n=d=2r+1$, $e=(1,\ldots,1)$, $X\in |2L-F|$;

\ind {\rm (d)} $r\geq 4$, $n=2r+1$, $d=2r$, $e=(1,\ldots,1)$, $X\in |2L-2F|$
or, equivalently, $X\simeq \PP^1 \times Q^{r-1}$ embedded Segre;

\ind {\rm (e)} $n=d=2r+2$, $e=(1,\ldots, 1,2)$, $X\in |2L-2F|$.
\end{thmno}

\begin{remno} {\rm (i) Except for case (i), all manifolds appearing in the theorem are rational.  

\hskip3pt (ii) All cases listed really occur.

(iii)  An inspection of the above list (or a direct argument)
shows that if we assume $d\leq r$, $X$ is either the Segre embedding of
$\PP^1 \times \PP^{r-1}$ in $\PP^{2r-1}$ or a Fano manifold with $b_2=1$.
In case $d\leq r-2$ all examples  I know of are complete intersections.

\hskip.5pt (iv) Manifolds from case (iv)~(b) up to (iv)~(e) in the theorem
are also Fano.}
\end{remno}

\section{Conventions and prerequisites}\label{S2}

We follow the customary notation in Algebraic Geometry (see e.g.~\cite{8}).
We denote by $X\subset \Pn _\mathbb{C}$ a complex projective
connected manifold. 
 We let $d$ be its degree and $r$ its dimension; $s=n-r$ is the
codimension of $X$ in $\Pn $. The irregularity of $X$ is
 by definition $q=:h^1(X, \co _X)$. $H$ will denote a
hyperplane section of $X\subset \Pn $. The sectional genus of $X$, denoted
$g$, is the genus of the curve $X\cap H_1\cap \cdots \cap H_{r-1}$, where
$H_1, \ldots, H_{r-1}$ are generic hyperplanes in $\Pn$. Adjunction
formula reads:
$$2g-2=(K+(r-1)H)H^{r-1},$$
where $K$ is a canonical divisor for $X$.

The $\Delta$-genus of $X$ is by definition
$$\Delta=d+r -h^0(X, \co _X (H))$$
and is a non-negative integer.

$X$ is said to be a {\it scroll} over the manifold $Y$ if $X\simeq \mathbb{
P}(\ce )$ for some vector bundle $\ce $ on $Y$, such that 
$\mathcal{O}_X(H)$ identifies to the tautological line bundle of $\PP (\mathcal{
E})$.

$X$ is said to be a {\it hyperquadric fibration} over the smooth curve
$C$ if there is a morphism $\pi: X\to C$ such that the fibres of $\pi$
are hyperquadrics with respect to the embedding induced by $\co _X
(H)$. It turns out that singular fibres of $\pi$ are ordinary cones (see
\cite{9}). In the sequel, we denote by $Q^r$ a hyperquadric of dimension~$r$.

The adjunction mapping of $X$, denoted below by $\varphi$, is the
rational map on $X$ associated to the linear system $|K+(r-1)H|$. See
e.g.\ \cite{3}, Chapters 9--11 for a complete study of its properties. 

We recall two results on the classification of manifolds of small
$\Delta$-genus. The first one is classical (see e.g.\ \cite{9}, Proposition 2.3).

\begin{thma}\label{TheoremA} The following are equivalent:

\hskip6pt{\rm (i)} $\Delta=0$;

\hskip3pt{\rm (ii)} $g=0$;

{\rm (iii)} $X$ is either $\PP ^r$, $H\in |\co_{\PP ^r}(1)|$ or the
hyperquadric $Q^r\subset \PP ^{r+1}$ or $\PP ^2$, $H\in
|\co_{\PP ^2}(2)|$ or a scroll over~$\PP ^1$.
\end{thma}

The next result is due to del Pezzo if $r=2$, to Fano and Iskovskih for
$r=3$ and to Fujita in general (see also \cite{9}, Proposition 2.4 for some other
characterisations).

\begin{thmb}\label{TheoremB} {\rm (Fujita, \cite{5}, \cite{6})} Assume that $r\geq 2$. The following
are equivalent:

\hskip3pt{\rm (i)} $\Delta=1$;

{\rm (ii)} $X$ is either a classical del Pezzo surface (anticanonical
embedding of either $\PP ^1\times \PP^1 $ or of the blowing-up of
$\PP ^2$ in at most six points) or , if $r\geq3$, one of the following:
 a hypercubic, a complete intersection of type $(2,2)$, a linear 
section of the Pl\" ucker embedding of the Grassmannian
of lines in $\PP ^4$, the Segre embedding of $\PP ^2\times \mathbb{
P}^2$, its hyperplane section (which is $\PP (T_{\PP ^2})$), the
Segre embedding of $\PP ^1\times \PP ^1\times \PP ^1$, the
scroll over $\PP ^2$, $\PP (\co _{\PP ^2} (1) \oplus 
\co _{\PP ^2}(2))$ (this is the blowing-up of $\PP ^3$ at a point), 
or the Veronese embedding~$v_2(\PP ^3)$.
\end{thmb}

Recall that $X$ is a Fano manifold if $-K$ is ample. We see that the
examples listed in Theorem B (which were called {\it classical del
Pezzo manifolds} in \cite{9}) are special Fano manifolds.  

\section{Proof of the theorem}\label{S3}

We begin with the following simple fact.
\begin{lem}\label{Lemma1} Let $C$ be a smooth projective curve of positive genus  and
let $\cl \in \Pic (C)$  with $\deg(\cl)>0$. Then we have
$h^0(\cl)\leq \deg(\cl)$.
\end{lem}
\begin{proof} If $\cl$ is special, we may apply Clifford's theorem.
If $\cl$ is non-special,
the result follows from the Riemann-Roch theorem. \end{proof}

\begin{prop}\label{Proposition1} Let $C$ be a smooth projective curve of positive genus
and let $\ce$ be an ample and spanned vector bundle on $C$. Then we have
$h^0(\ce)\leq\deg(\ce)$.
\end{prop}
\begin{proof} We proceed by induction on $e=:\rank(\ce)$. When $e=1$, we may apply
Lemma~\ref{Lemma1}. Assume now $e\geq 2$. As $\ce$ is ample and spanned, it follows
that $h^0(\ce)>e$. So, for $p\in C$, we may find a non-zero section
$s\in H^0(C, \ce(-p))$. $s$ induces an exact sequence:
$$0 \longrightarrow \cl \longrightarrow\ce \longrightarrow \ce'
\longrightarrow 0,$$
where $\cl \in \Pic(C)$, $\deg (\cl)=:l>0$, and $\ce'$ is ample, spanned
and of rank $e-1$. We have
$$\deg(\ce)-l=\deg(\ce')\geq h^0(\ce')\geq h^0(\ce)-h^0(\cl)$$
by the induction hypothesis and the cohomology sequence of the above
exact sequence. Applying once again Lemma~\ref{Lemma1} we get $\deg (\ce)\geq
h^0(\ce)$.\end{proof}
\begin{cor}\label{Corollary1} Let $X\subset \Pn$ be a scroll over a smooth curve
$C$. Assume that $X$ is non-degenerate and $d\leq n$. Then $C\simeq
\PP^1$.
\end{cor}
\begin{proof} Let $X\simeq \PP (\ce)$. If $g(C)>0$, by Proposition~\ref{Proposition1} we get
$$n+1\leq h^0(X, \co_X(H))=h^0(C, \ce)\leq \deg (\ce)=d,$$ a
contradiction. \end{proof}
\begin{lem}\label{Lemma2} Let $X\subset \PP^{n=r+s}$ be smooth connected 
non-degenerate with $d\leq n$. Assume moreover that $r\leq
s+1$. Then we
have:

\hskip3pt{\rm (i)} $g\leq r-1$; and

{\rm (ii)} $d\geq 2g+1$.
\end{lem}
\begin{proof} (i) Let $C\subset \PP^{s+1}$ be the curve-section of $X$. If
$H_C$ is special, by Clifford's theorem we get
$$s+2\leq h^0(C,\co_C(H_C))\leq \frac{d}{2} +1\leq \frac{r+s}{2}+1,$$ giving
$r\geq s+2$. This is a contradiction. So $H_C$ is non-special and by
Riemann-Roch we get
$$s+2\leq h^0(C, \co_C(H_C))=d+1-g\leq r+s+1-g,$$
hence $g\leq r-1$.

(ii) Assume that $d\leq 2g$. We get by (i)
$$r+1\leq s+2\leq h^0(C, \co_C(H_C))=d+1-g\leq g+1\leq r,$$
which is absurd. \end{proof}
\begin{prop}\label{Proposition2} Let $X\subset \Pn$ be smooth connected 
non-degenerate and linearly normal with $d\leq n$. Assume that the
adjunction mapping $\varphi =\varphi_{|K+(r-1)H|}$ makes $X$ a scroll
over a smooth surface $S$. Then  $S\simeq \PP^2$ and $X$ is
one of the following:

$r=4$, $d=10$, $X\simeq \PP(T_{\PP^2}\oplus\co_{\PP^2}(1))$, or

$r=4$, $d=11$, $X\simeq \PP(\co_{\PP^2}(1)\oplus \co_{\PP^2}(1)\oplus
\co_{\PP^2}(2))$, or
$r=5$, $d=10$, $X\simeq \PP(\co_{\PP^2}^{\oplus 4}(1))$, i.e. $X$ is the
Segre embedding of $\PP^2\times \PP^3$.
\end{prop}
\begin{proof} Let $S'$ be the smooth surface
$X\cap H_1\cap\cdots\cap H_{r-2}$, where $H_i$ are generic hyperplanes
in $\PP^n$. We first remark that the geometric genus of $S'$ is zero. This
follows from Lemma~\ref{Lemma2}~(ii) and the adjunction formula for $H_{S'}$. The
standard exact sequences
\begin{align*}0 \longrightarrow\co_X(K+(r-2)H)& \longrightarrow\co_X(K+(r-1)H)
\\ &\longrightarrow \co_H(K_H+(r-2)H_H)\longrightarrow 0\end{align*}
together with Lemma 1.1 from \cite{9} show that, in our case,
$h^0(X,\co_X(K+(r-1)H))=g-q$. So, we have $\varphi:X\to S\subset \PP^{g-
q-1}$.
Let $H_S$ be the generic hyperplane section of $S\subset\PP^{g-q-1}$ and
let $Y=:\varphi^{-1}(H_S)$. Note that $Y$ is a scroll of dimension $r-1$
over the curve $H_S$; if we let $d_Y$ be its degree, we get $d_Y= (K+(r-
1)H)H^{r-1}=2g-2$ by adjunction formula. Let $m$ be the dimension of the
projective space spanned by $Y$ inside $\PP^n$ (denoted below by 
$\langle Y \rangle$).
 By Barth's
theorem (see \cite{2}) we must have $m\geq 2(r-1)-1$. We get, using Lemma~\ref{Lemma2}~(i)
$$m\geq 2r-3\geq 2(r-2)\geq 2(g-1)=d_Y.$$
So, by Corollary~\ref{Corollary1}, it follows that $H_S\simeq \PP^1$. The 
two-dimensional version of Theorem~A shows that $q=0$ and one of the
following holds:

\vskip6pt
1. $S=\PP^2$, $g=\Delta=3$;

2. $S$ is a scroll over $\PP^1$;

3. $S$ is the Veronese embedding $v_2(\PP^2)$, $g=6$.

\vskip6pt
\noindent Recalling the definition of $\Delta=d+r-h^0(X,\co_X(H))$, we get
$$n+r\geq d+r\geq n+1+\Delta,$$
giving $r\geq \Delta +1$.
Now, if we are in case 1, by Proposition 4.7 from \cite{9}, it follows that we
have the following possibilities for $X$:

$r=4$, $d=9,10$ or $11$;

$r=5$, $d=10$, $X$ is the Segre embedding of $\PP^2\times \PP^3$.

\noindent Assume that $r=4$, so $X\simeq \PP (\ce)$ for some very ample
vector bundle of rank three over $\PP^2$. If $\ell$ is a line in
$\PP^2$, it follows that $\ce|_\ell$ has degree $4$ and is very ample.
So, $\ce|_\ell\simeq \co_{\PP^1}(1)\oplus\co
_{\PP^1}(1)\oplus\co_{\PP^1}(2)$, i.e. $\ce$ is uniform. One may use the
classification from \cite{4}; we find  that case $d=9$ is not
possible, while for $d=10$ we get $\ce \simeq
T_{\PP^2}\oplus\co_{\PP^2}(1)$ (equivalently $X$ is the hyperplane
section of the Segre embedding of $\PP^2\times \PP^3$) and for $d=11$ we
get $\ce \simeq \co_{\PP^2}(1)\oplus \co_{\PP^2}(1)\oplus \co
_{\PP^2}(2)$ (this is the blowing-up of $\PP^4$ with center a line).

To finish the proof we only have  to show that cases 2 and 3 cannot occur. 
We use the
notation from \cite{8}, Chapter V, Section 2. If we are in case 2, we have
$S\simeq \mathbb{F}_e$, $H_S=C_0+bF$ with $b>e\geq 0$.

We look at the $(r-1)$-dimensional rational scrolls $Y_0=\varphi^{-
1}(C_0)$ and $Y_1=\varphi^{-1}(F)$. If we put $d_i =\deg (Y_i)$ for
$i=0,1$, we get $d_i\geq r-1$. Indeed, by Barth's theorem (\cite{2}), if
$m_i=\dim\langle Y_i \rangle$, we get $m_i\geq 2(r-1)-1$; moreover,
since $\Delta (Y_i, \co_{Y_i}(H))=0$, we deduce
$$d_i+r-1=h^0(Y_i, \co_{Y_i}(H))\geq m_i+1\geq2(r-1),$$
i.e.\ $d_i \geq r-1$. So, we find
$$2g-2=\deg(Y)=d_0+bd_1\geq d_0+d_1 \geq 2(r-1),$$
 contradicting part (i) of  Lemma~\ref{Lemma2}. Case 3 is ruled out  by a similar argument.\end{proof}

Next we need a general lemma concerning 
the geometry of hyperquadric fibrations (see also \cite{10}, 6.2).

\begin{lem}\label{Lemma3} Assume that the adjunction mapping 
$\varphi : X\to C\subset \PP^m$ makes 
$X$ a hyperquadric fibration over the smooth curve $C$.
Then $m=g-q-1$ and $q$ coincides with  the genus of $C$.
 Moreover, if we let $\ce=: \varphi_*\co_X(H)$, $\ce$ is a spanned
vector bundle of rank $r+1$ over $C$. Denote by $\pi: \PP(\ce)\to C$
the projection and by $L$ the tautological divisor on $\PP(\ce)$.
Then $X$ is embedded in $\PP (\ce)$ such that $L|_X=H$ and $X\in |2L+\pi^*B|$ for
some divisor $B$ on $C$.
Finally, if $a =:\deg(\ce)$ and $b=:\deg (B)$, the following formulae
hold
$$a= 1-g+2(q-1)+d\quad \hbox{and}\quad b=2(g-1)-4(q-1)-d.$$
\end{lem}
\begin{proof}
From Lemma 1.1 in \cite{9} and the standard exact sequences
\begin{align*}0\longrightarrow \co_X(K+(r-2)H)&\longrightarrow
 \co_X(K+(r-1)H)\\ &\longrightarrow 
\co_H(K_H+(r-2)H_H)\longrightarrow  0\end{align*}
\vskip3pt

\noindent it follows in our case that $h^0(X,\co_X(K+(r-1)H))=g-q$. We have, for
any $c\in C$, $H^0(X_c, \co_{X_c}(H))=r+1$ and $H^1(X_c,
\co_{X_c}(H))=0$, so the existence of $\ce$ follows from the base-change
theorem. Moreover, the canonical diagram
\vskip3pt
$$\begin{matrix}H^0(C,\ce)&{\buildrel \sim\over \longrightarrow}& H^0(X,
\co_X(H))\cr
\noalign{\vskip4pt}
{\scriptstyle{\rm ev}} \Big\downarrow && 
\Big
\downarrow{\scriptstyle\,{\rm res}}
\cr
\noalign{\vskip4pt}\ce_c&
 {\buildrel \sim\over \longrightarrow} & H^0(X_c, \co_{X_c}(H))\end{matrix}
$$
\vskip8pt
\noindent shows that $\ce$ is spanned by global sections iff the restriction map
res is surjective for any $c\in C$. This holds true since $X_c$ is a
hyperquadric, hence linearly normal in $\PP^r=\langle X_c\rangle$.

Consider also the canonical induced diagram
$$\begin{matrix}X &\subset&\PP(\ce)\cr\noalign{\vskip4pt}
\,\,\,\,{\scriptstyle\varphi}\searrow&&\swarrow {\scriptstyle\pi}\,\,\,\,\,\cr
\noalign{\vskip4pt}
&C&\end{matrix}$$
and write $X\sim 2L+\pi^*B$, for some $B\in {\rm Div}(C)$. Let $H_C$ be
the hyperplane section of $C\subset \PP^{g-q-1}$. We find
$$ \varphi^*(H_C)=K+(r-1)H=(K_{\PP(\ce)}+X+(r-
1)L)|_X=\varphi^*(K_C+\det\ce +B).$$
By taking degrees, we get $g-1=2(q-1)+a+b$. Moreover, $a=(L^{r+1})$, so
$d=(L^r\cdot X)=2a+b$.
The two formulae follow. \end{proof}
\begin{lem}\label{Lemma4} Let $X\subset \Pn$ be smooth connected non-degenerate
with $d\leq n$. Assume that the adjunction mapping $\varphi:X\to C$
makes $X$ a hyperquadric fibration over the smooth curve $C$.
Then $C\simeq\PP^1$.
\end{lem}
\begin{proof} Assume that $q=g(C)>0$. By Lemma~\ref{Lemma2}~(ii), $d\geq 2g +1$. So, by
Lemma~\ref{Lemma3}, we have $b=2(g-1)-d-4(q-1)<0$.

We show first that $\ce$ is ample. As $\ce$ is spanned,
$\co_{\PP(\ce)}(L)$ is spanned. So, if $L$ is not ample, there is a
curve $D\subset \PP(\ce)$ such that $(L\cdot D)=0$. It follows that
$(X\cdot D)=(2L+\pi^* B)D= \alpha b$ for some $\alpha>0$.
As $b<0$, we deduce that $(X\cdot D)<0$, so $D\subset X$. But $L|_X=H$,
so $(D\cdot L)>0$ which is a contradiction. So $\ce$ is ample.

Let now $S\subset X$ be the surface-section of $X$, i.e. $S=X\cap H_1
\cap\cdots \cap H_{r-2}$, where $H_i$ are generic hyperplanes in $\Pn$.
We have $(H_S+K_S)^2 =0$, giving $d+2(H_S\cdot K_S) +(K_S)^2=0$.
Adjunction formula yields $(H_S\cdot K_S)= 2g-2-d$; moreover,
$(K_S)^2\leq 8(1-q)$, since $S$ is birationally ruled. We deduce, using
also Lemma~\ref{Lemma2}~(ii)
$$4(g-1)\geq d+8(q-1)\geq 2g+1 +8(q-1).$$
So we get $4q\leq g+1.$ By Lemma~\ref{Lemma3}, $a=1-g+2(q-1)+d$ and we find
$a\leq d-2q$. Now, since $\ce$ is ample and spanned, we may apply 
Proposition~\ref{Proposition1} to find
$$a=\deg (\ce)\geq h^0(C,\ce)=h^0(X, \co_X(H))\geq n+1.$$
Putting things together, we get
 $$n+1\leq a\leq d-2q\leq n-2q.$$
This is a contradiction, so $q=0$. \end{proof}

We shall also need the proposition below which might  have an interest in itself.
\begin{prop}\label{Proposition3} Let $X\subset \Pn$ be smooth, connected, 
non-degenerate and linearly normal. Assume that the adjunction 
mapping $\varphi:X\to C$ makes $X$ a 
hyperquadric fibration over $C\simeq \PP^1$. Assume moreover, that
$d\geq 2g+2$ and $r\geq g+1$. Then, in the notation of Lemma~{\rm \ref{Lemma3}} and
denoting by $e=(e_0, \ldots, e_r)$ the splitting type of $\ce$ and by
$F$ a fibre of the projection $\PP(\ce)\to \PP^1$, we have one of the
following:

{\rm (a)} $r=s$, $d=2r$, $e=(1,\ldots, 1,0)$, $X\in |2L|$;

{\rm (b)} $r=s-1$, $d=2r+1$, $e=(1,\ldots, 1)$, $X\in |2L-F|$;

{\rm (c)} $r=s-1$, $d=2r$, $e=(1,\ldots, 1)$, $X\in |2L-2F|$ or,
 equivalently, $X\simeq \PP^1\times Q^{r-1}$ embedded Segre;

{\rm (d)} $r=s-2$, $d=2r+2$, $e=(1, \ldots, 1, 2)$, $X\in |2L-2F|$;

{\rm (e)} $r=3$, $X\simeq \PP^1\times \mathbb{F}_1$, embedded Segre, where 
$\mathbb{F}_1$ is embedded in $\PP^4$ as a rational scroll of degree $3$.

\noindent Moreover, all these cases do exist.
\end{prop}
\begin{proof} We first remark that $g\geq 2$ (see \cite{9}), so $r\geq 3$. Let 
$Q$ denote a fibre of $\varphi$. We have $(H-Q)H^{r-1}=d-2$.  The
standard exact sequence
$$0\longrightarrow \co_X(-Q)\longrightarrow \co_X(H-Q) 
\longrightarrow \co_H(H-Q)\longrightarrow 0$$
and the fact that $H^1(X,\co_X(-Q))=0$ allow one to prove by induction
on $r$ that $|H-Q|$ is base-points free. Note that on the curve-section
of $X$, the degree of the restriction of $|H-Q|$ is $\geq 2g$, 
so it is base-points free. Moreover, $|H-Q|$
is not composed with a pencil, since $r\geq 3$. So, by Bertini's
theorem, there is a smooth member $X'\in |H-Q|$. We let
\begin{alignat*}{2}
&H' =H|_{X'}, &\qquad& K'=K_{X'},\\ &r'=\dim(X')=r-1,&
&\varphi'=\varphi _{|K' + (r'-
1)H'|},\\
&d'=\deg(X')=d-2, &&g'=g(H'), \\&s'=h^0(X', \co_{X'}(H))-1-
r'.\end{alignat*}
One finds easily $g'=g-1$, $s'=s-1$ and $\varphi'$ identifies to
$\varphi|_{X'}$. The statement of the proposition is proved by induction
on $r$ (note that we still have $d'\geq 2g'+2$ and $r'\geq g'+1$).
Assume first that $g\geq 3$. Since $r\geq g+1$, for $r=4$ we get $g=3$
and we may use the classification from Theorem 4.3 in \cite{9}. For $r\geq4$
we find inductively the following possible values for the numerical
invariants:
\vskip6pt

(a) $r=s$, $d=2r$, $g=r-1$;

(b) $r=s-1$, $d=2r+1$, $g=r-1$;

(c) $r=s-1$, $d=2r$, $g=r-2$;

(d) $r=s-2$, $d=2r+2$, $g=r-1$.
\vskip6pt

\indent It remains to analyse the case $g=2$, where one may use the
classification theorem 3.4 in \cite{9}. This leads to only one new case, which is
(e).

Next we investigate the structure of $\ce $ in each case.

First we have that $\ce$ is non-special (since it is spanned by  Lemma~\ref{Lemma3}).
So Riemann-Roch theorem gives
$$r+s+1=h^0(\ce)=a +r+1,$$
hence $a=s$. Now, in case (a), we remark that $|H-2Q|=\emptyset$, since
$(H-Q)^{r-1}\cdot(H-2Q)=d-2r-2<0$.
By Lemma~\ref{Lemma3}, $b=0$, so $X\in |2L|$.

The exact sequence
$$0\longrightarrow \co_{\PP(\ce)}(-L-2F) \longrightarrow
\co_{\PP(\ce)}(L-2F)
\longrightarrow \co_X(H-2Q)\longrightarrow 0$$
shows that $h^0(\ce(-2))=0$; as $\ce$ is spanned and $a=r$, the
splitting-type of $\ce$ must be $(1,\ldots,1,0)$. The existence follows
by the same type of argument as in the proof of Proposition 3 from \cite{12}.
The other cases are similar and simpler. For instance, in case (b) one 
gets as above $h^0(\ce(-2))=0$, $a=r+1$ and $b=-1$. So $e=(1,\ldots,1)$, 
$\ce$ is very ample and the existence follows now easily. \end{proof}

\begin{prop}\label{Proposition4} Let $X\subset \Pn$ be smooth, connected, 
non-degenerate and linearly normal, with $d \leq n$. Assume that the 
adjunction mapping makes $X$ a hyperquadric fibration over a smooth curve
$C$. Then $X$ is as in case {\rm (ii)~(b)} or case {\rm (iv)} of the main
theorem.
\end{prop}
\begin{proof}
By Lemma~\ref{Lemma4} $C\simeq \PP^1$. We have $d\geq 2g+1$ and $g\leq r-1$ by
Lemma~\ref{Lemma2}. If $d\geq 2g+2$, we may apply Proposition~\ref{Proposition3}, thus leading to
cases (ii)~(b) and (iv)~(b) up to (iv)~(e) of the main theorem.
So, assume that $d=2g+1$. As in the proof of Proposition~\ref{Proposition3} we deduce
that $a=s$. By Lemma~\ref{Lemma3} we get $a=g$, $b=1$.
 It follows $s=g\leq r-1$. Barth's theorem (\cite{2}) ensures that $s\geq r-1$,
so we must have $s=r-1$. We obtain 
$$g=r-1, \quad d= 2r-1,\quad a=r-1.$$
As in the proof of Proposition~\ref{Proposition3}, we have $|H-2Q|=\emptyset$, so
$h^0(\ce(-2))=0$.
It follows that the splitting type of $\ce$ is $(1, \ldots, 1, 0, 0)$,
so we are in case (iv)~(a) of the main theorem. The existence follows
from Proposition~3~in~\cite{12}.\end{proof}

We are now ready for the proof  of our theorem.

Assume first that $r\leq s+1$. We have
$$\Delta=d+r-h^0(X_, \co_X (H))\leq n+r-n-1=r-1.$$
If $\Delta =0$, by Theorem A we get either case (iii) of the main
theorem or some special examples of case (i). Similarly if $\Delta =1$,
by Theorem B we get either case (ii)~(a) or some special examples of 
case (i). So, assume $\Delta \geq 2$, hence $r\geq 3$, from now on. 
If $r=3$, it follows $\Delta=2$, $s\geq 2$ and $\varphi :X\to \PP^1$ is a 
hyperquadric fibration by \cite{9}, Theorem~3.12 and Corollary~3.3. If $r=4$, 
we get $\Delta =2$ or $3$, $s\geq 3$, so $\varphi$ is either a
hyperquadric fibration over a rational curve or a scroll over $\PP^2$
(see \cite{9}, Theorems~3.12, 4.8
 and 4.2).
Since $d\leq n$, it follows that $d\leq r+s\leq 2s+1$. So, using the
general properties of the adjunction mapping (see e.g. \cite{3}, Chapters 9-11, in
particular Theorem 11.2.4)
 and the above analysis for $r\leq 4$, it follows
from Theorem I in \cite{11} that one of the following holds:
\vskip6pt

(1)
$X$ is a scroll over a (smooth) curve $C$;

(2) $\varphi$ makes $X$ a scroll over a smooth surface;

(3) $\varphi$ makes $X$ a hyperquadric fibration over a smooth curve.
\vskip6pt

In case (1), from Corollary~\ref{Corollary1}, we get $C\simeq \PP^1$, 
so $\Delta =0$. In case (2), by Proposition~\ref{Proposition2} we reach 
case (ii)~(c). If we are in case (3),
 by Proposition~\ref{Proposition4} we get case (ii)~(b) or case (iv). Assume now that $r\geq s+2$. 
By Barth's theorem (\cite{2}) it follows that 
$\Pic (X)\simeq \mathbb{Z}$, generated by the class 
of $\co_X(H)$. We show that $X$ is Fano, so we are in case (i) and 
the main theorem is 
completely proved. As we have $\Pic (X)\simeq  \mathbb{Z}$, to prove that
$X$ is Fano it is enough to see that the geometric genus of $X$, denoted
$p_g$, is zero. Here we make use of a theorem of Harris (see \cite{7}),
generalising Castelnuovo's bound for the genus of a curve to
arbitrary dimension.

It states that 
$${\rm p_g}\leq \binom{M}{r+1}s + \binom{M}{r} \varepsilon,$$
where $M=[(d-1)/s]$ and $\varepsilon =d-1-Ms$.

If $s=1$ we find ${\rm p_g}=0$  by an obvious direct computation.
If $s\geq 2$ and $r\geq 2$ we get $r+s-1<rs$; our hypothesis $d\leq r+s$ 
gives $d-1< rs$, or $M<r$. So ${\rm p_g}=0$.\qed

\section*{Acknowledgements}  

I thank Francesco Russo for 
sending me the interesting preprint (\cite{1}) and especially for  
explicitely asking me if I know examples of (non-degenerate) irregular
$r$-dimensional manifolds in $\PP^{2r+1}$, of degree $\leq 2r+1$. In fact his 
question was the starting point of the present work.

\end{document}